\documentclass[11pt]{article}
\usepackage{amsmath}
\usepackage{amsfonts}
\usepackage{amsthm}
\usepackage{verbatim}
\def\N{\mathbb{N}}
\def\F{\mathbb{F}}

\def\H{\mathcal{H}}

\begin{document}
\centerline{\bf{ HYPERQUADRATIC CONTINUED FRACTIONS }}
\centerline{\bf{ OVER A FINITE FIELD OF ODD CHARACTERISTIC }}
\centerline{\bf{ WITH PARTIAL QUOTIENTS OF DEGREE 1}}
\vskip 0.5 cm
\centerline{\bf{by}}
\centerline{\bf{A. Lasjaunias and J.-Y. Yao}}
\vskip 0.5 cm {\bf{Abstract.}} In 1986, some examples of algebraic, and
nonquadratic, power series over a finite prime field, having a continued fraction expansion with partial quotients all of degree one, were discovered by
W. Mills and D. Robbins. In this note we show how these few examples are
included in a very large family of continued fractions for certain
algebraic power series over an arbitrary finite field of odd characteristic. 
\vskip 0.5 cm
\noindent{\bf{1. Introduction and results}}
\par The subject of this note belongs to the theory of continued fractions
in power series fields. For a general account on this matter, the
reader can consult W. Schmidt's article [14]. For a wider survey on diophantine approximation in the function fields case and full references, the reader may also consult D. Thakur's book [15, Chap. 9]. Let us recall that the pioneer work on the matter treated here, i.e. algebraic continued fractions in power series fields over a finite field, is due to L. Baum and M. Sweet [2]. 
\par Let $p$ be a prime number, $q=p^s$ with $s\geq 1$, and let $\F_q$ be the finite field with $q$ elements. We let $\F_q[T]$, $\F_q(T)$
 and $\F(q)$  respectively denote,  the ring of polynomials, the field of rational functions and the field of power series in $1/T$
over $\F_q$, where $T$ is a formal indeterminate. These fields are
valuated by the ultrametric absolute value (and its extension) introduced on $\F_q(T)$ by
$\vert P/Q\vert=\vert T\vert^{\deg(P)-\deg(Q)}$, where $\vert T\vert>1$ is a fixed real number. We recall that each
irrational (rational) element $\alpha$ of $\F(q)$ can be expanded as an infinite
(finite) continued fraction. This will be denoted
$\alpha=[a_1,a_2,\dots,a_n,\dots]$ where the $a_i\in \F_q[T]$, with
$\deg(a_i)>0$ for $i>1$, are the
partial quotients and the tail $\alpha_i=[a_i,a_{i+1},\dots]\in\F(q)$ is the complete quotient. We shall be concerned with infinite continued fractions in $\F(q)$ which are algebraic over $\F_q(T)$. 
\par Regarding diophantine approximation and continued fractions, a particular subset of elements in $\F(q)$, algebraic over $\F_q(T)$, must be considered. Let $r=p^t$ with $t\geq 0$, we denote by $\H(r,q)$ the subset of irrational $\alpha$ belonging to $\F(q)$ and satisfying an algebraic equation of the particular form $A\alpha^{r+1}+B\alpha^r+C\alpha+D=0$, where $A,B,C$ and $D$ belong to $\F_q[T]$. Note that $\H(1,q)$ is simply the set of quadratic irrational elements in $\F(q)$. The union of the subsets $\H(p^t,q)$, for $t\geq 0$, denoted by $\H(q)$, is the set of hyperquadratic power series. For more details and references, the reader may see the introduction of [4]. Even though it contains algebraic elements of arbitrary large degree, this subset $\H(q)$ should be regarded as an analogue, in the formal case, of the subset of quadratic numbers, in the real case. An old and famous theorem, due to Lagrange, gives a characterization of quadratic real numbers as ultimately periodic continued fractions. It is an open problem to know whether another characterization, as particular continued fractions, would be possible for hyperquadratic power series.
\par The origin of this work is certainly due to a famous example of a cubic power series over $\F_2$, having partial quotients of bounded degrees (1 or 2), introduced in [2]. In a second article [3], Baum and Sweet could characterize all power series in $\F(2)$ having all partial quotients of degree 1 and, among them, those which are algebraic. Underlining the singularity of this context, in [9, p. 5], a different approach could allow to rediscover these particular power series in $\F(2)$. Also in characteristic 2, other algebraic power series over a finite extension of $\F_2$, having all partial quotients of degree 1, were presented (see for instance [10, p. 280]). The case of even characteristic appears to be singular for different reasons. In this note we only consider the case of odd characteristic. Our aim is to show the existence of hyperquadratic continued fractions, in $\F(q)$ for $p\neq 2$, with all partial quotients of degree 1. The first examples, in $\F(p)$, were presented by Mills and Robbins [12]. 
\par Before developing the background of the work presented in this article, we first give an example of such algebraic continued fractions with the purpose of illustrating the subject discussed here. The following result is derived from an elementary and particular case of the theorem which is stated at the end of this section. 
\vskip 0.5 cm
\noindent {\bf{Example.}}{\emph{  Let $p$ be an odd prime number. Let $\epsilon \neq 0,1$ in $\F_p$. Let us consider the algebraic equation, with coefficients in $\F_p[T]$:
$$X^{p+1}-TX^p+\epsilon T((T^2-1)^{(p-1)/2}-T^{p-1})X+\epsilon (T^{p+1}-(T^2-1)^{(p-1)/2}(T^2+\epsilon-1))=0.$$
This equation has a unique root $\alpha$ in $\F(p)$, with $\vert \alpha\vert \geq
\vert T\vert$,  which can be expanded as the following infinite continued fraction 
$$\alpha=[T,(\epsilon (\epsilon-1))^{-1}T,(2\epsilon T,-2\epsilon^{-1}T)^{(p-1)/2},\dots,(\epsilon (\epsilon-1))^{u_m}T,(2v_mT,-2v_m^{-1}T)^{(p^m-1)/2},\dots \dots],$$
where $(a,b)^k$ denotes the finite sequence $a,b,a,\dots,b$ of length $2k$, the pair $a,b$ being repeated $k$ times, with $u_m=-1$ if $m$ is odd and $u_m=0$ if $m$ is even, while $v_m=\epsilon$ if $m$ is odd and $v_m=(\epsilon-1)^{-1}$ if $m$ is even.}} 
\par To explain the existence of such continued fractions, our method is based on the following statement, proved by the first author [6, pp. 332-333].
\par {\it{Given an integer $l\geq 1$, a $l$-tuple $(a_1,a_2,\dots,a_l) \in (\F_q[T])^l$, with $\deg(a_i)\geq 1$ for $1\leq i\leq l$,
 and a pair $(P,Q)\in (\F_q[T])^2$ with $deg(Q)<deg(P)<r$, there exists a unique infinite continued fraction $\alpha \in \F(q)$ 
defined by
$$\alpha=[a_1,\dots,a_l,\alpha_{l+1}] \quad \text{ and } \quad \alpha^r=P\alpha_{l+1}+Q. \eqno{(*)}$$}}
\par Note that in the degenerated case, $r=1$, consequently $\deg(P)=0$ and $Q=0$, we simply have $\alpha=\epsilon \alpha_{l+1}$, where $\epsilon \in \F_q^*$. This implies the (pure) periodicity of the continued fraction, with a period of length multiple of $l$. In general, the continued fraction $\alpha$ is algebraic over $\F_q(T)$ of degree $d$,
with $1<d\leq r+1$. Indeed, from the continued fraction algorithm,  we know that
there is a linear fractional transformation $f_l$, having coefficients in $\F_q[T]$, built from the $l$ first partial
quotients, such that
 $\alpha=f_l(\alpha_{l+1})$ (see the end of Section 2). Consequently, by $(*)$  we have
 $\alpha=f_l((\alpha^r-Q)/P)=f(\alpha^r)$ where $f$ is a linear fractional
 transformation with integer (polynomial) coefficients. Hence,
 $\alpha$ is solution of the following algebraic equation of degree $r+1$: 
$$y_lX^{r+1}-x_lX^r+(Py_{l-1}-Qy_l)X+Qx_l-Px_{l-1}=0,\eqno{(**)}$$
where the polynomials $x_l,x_{l-1},y_l$ and $y_{l-1}$ are the
continuants built from the $l$ first partial quotients (see the end of Section 2). Thus, $\alpha$ is hyperquadratic. Moreover, it is also true that $\alpha$ is the unique root in $\F(q)$, with $\vert \alpha\vert \geq
\vert T\vert$, of equation $(**)$.
\par In this note, we shall consider continued fractions in $\F(q)$
defined by $(*)$, for a particular choice of the polynomials
$(a_1,a_2,\dots,a_l,P,Q)$. Here we consider $p>2$, $q=p^s$ and $r=p^t$ as above. In the sequel $a$ is given in $\F_q^*$. We
consider the following pair of polynomials in $\F_q[T]$:
$$P_a(T)=(T^2+a)^{(r-1)/2} \quad \text{ and } \quad Q_a(T)=a^{-1}(TP_a(T)-T^r).$$
We have $\deg(P_a)=r-1>\deg(Q_a)=r-2$. For an integer $l\geq 1$, we denote by $\mathcal{E}(r,l,a,q)$ the
subset of infinite continued fraction expansions $\alpha
\in \F(q)$ defined by
$$\alpha=[a_1,\dots,a_l,\alpha_{l+1}] \quad \text{ and}\quad \alpha^r=\epsilon_1P_a\alpha_{l+1}+\epsilon_2Q_a,$$ 
where $a_i=\lambda_iT+\mu_i$, $\lambda_i \in
\F_q^*$, $\mu_i \in \F_q$, for $1\leq i\leq l$ and $(\epsilon_1,\epsilon_2)\in (\F_q^*)^2$
are given. Note that in the extremal case, $r=1$, the pair of polynomials would be $P_a=1$ and $Q_a=0$ and $\mathcal{E}(1,l,a,q)$ would be a subset of quadratic power series, corresponding to purely periodic continued fractions. In the sequel we assume $r>1$. We observe that $\alpha$ in $\mathcal{E}(r,l,a,q)$ is
defined by the $(2l+2)$-tuple
$(\lambda_1,\dots,\lambda_l,\mu_1,\dots,\mu_l,\epsilon_1,\epsilon_2)$. Consequently $\mathcal{E}(r,l,a,q)$ has $q^l(q-1)^{l+2}$ elements.
\newline Our aim is to show that, under a particular choice of the $l-$tuple
$(a_1,\dots,a_l)\in (\F_q[T])^l$ and of the pair
$(\epsilon_1,\epsilon_2)\in (\F_q^*)^2$, the element $\alpha$ defined
as above will satisfy $\deg(a_n)=1$, for all the partial quotients $a_n$ in its continued fraction expansion. These particular
expansions are said perfect and form a subset of $\mathcal{E}(r,l,a,q)$, which will be denoted by $\mathcal{E^*}(r,l,a,q)$.
\newline A particular and simpler case of this situation can be considered. Let us denote by $\mathcal{E}_0(r,l,a,q)$ the subset of $\mathcal{E}(r,l,a,q)$
where $a_i=\lambda_iT$, for $1\leq i\leq l$, and also
$\mathcal{E}^*_0(r,l,a,q)=\mathcal{E^*}(r,l,a,q)\bigcap
\mathcal{E}_0(r,l,a,q)$. Considering the algebraic equation which they
satisfy, it can be observed that the continued fractions
belonging to $\mathcal{E}_0(r,l,a,q)$ are odd functions of $T$ and
therefore the partial quotients must be odd polynomials of the
indeterminate $T$. Consequently, the elements of $\mathcal{E}^*_0(r,l,a,q)$ have all partial
quotients of the form  $a_n=\lambda_nT$, for $n\geq 1$, where $\lambda_n
\in \F_q^*$. It can be observed that the example introduced above actually belongs to $\mathcal{E}^*_0(p,1,-1,p)$, and it is defined by the triple $(\lambda_1,\epsilon_1,\epsilon_2)=(1,\epsilon(\epsilon-1),\epsilon)$.
\par In Mills and Robbins article [12], several examples of algebraic
continued fractions are
presented, some of them with all partial quotients of degree 1. The first examples, [12, pp. 400-401], belong to
$\mathcal{E}^*_0(p,2,4,p)$, for all primes $p\geq 5$ (see [6, p. 332]). Also in [12, pp. 401-402], we have an example belonging to $\mathcal{E^*}(3,7,1,3)$. In this last case, the partial quotients are not linear. Inspired by this
example and using a new approach, in an earlier work [8], the first author could present a particular family of such hyperquadratic continued fractions, all in $\mathcal{E^*}(3,l,1,3)$, for all $l\geq 3$. In a joint work with J-J. Ruch [10], a generalization of the approach introduced in [8] was developed, for all characteristics; however, this led to unsolved questions. 
\par Yet, other examples of algebraic
continued fractions, with partial quotients of unbounded degrees, were also presented by Mills and Robbins [12]. One could observe that all these continued fractions are generated as indicated above, in connection with a polynomial of the form $(T^2+a)^k$, for different values of an integer $k$. Following this, in a larger context than the one we consider here (i.e. for different values of $k$, see [6] and [7] ), the first author could develop a method showing the link between all these algebraic continued fractions. In a particular case, this method can be used to describe continued fractions in $\mathcal{E}^*_0(r,l,a,q)$, including the previously mentioned examples [12, pp. 400-401]. However this method (introduced in [6] and developped in [7]) concerned only elements in $\mathcal{E}_0(r,l,a,q)$. In the simplest case, $p=q=r=3$, in a recent joint work with D. Gomez [5], a modification of this approach has allowed to obtain a large extension of the results presented in [8]. The aim of the present note is to give a full description of these particular algebraic continued fractions for all $p>2$, $r$ and $q$.
\par Before stating our result, it is pertinent, just for the sake of completeness, to recall what is already known in this area. Indeed, the following could be proved [7, p. 256]:
\vskip 0.5 cm
\noindent
 {\it{If $\alpha \in \mathcal{E}_0(r,l,-1,q)$ and if we have 
  $(C_0)$ : $[\lambda_1,\lambda_2,\dots,\lambda_l+\epsilon_1/\epsilon_2]^r=\epsilon_2$,
  then   $\alpha \in \mathcal{E}^*_0(r,l,-1,q)$.}} 
\vskip 0.5 cm
\noindent This condition must be understood as a set of several conditions
implying the existence of the square bracket on the left. As an illustration, for $l=1$, $(C_0)$ is simply
$(\lambda_1+\epsilon_1/\epsilon_2)^r=\epsilon_2$. There are
$q-1$ choices for $\epsilon_2$ and, for each one, $q-2$ choices for
$\lambda_1$, since $\lambda_1^r\neq 0,\epsilon_2$, while $\epsilon_1$ is fixed by
$\epsilon_1^r=(\epsilon_2-\lambda_1^r)\epsilon_2^r$. More generally, we can observe that there are $(q-1)^{l+2}$ elements in
$\mathcal{E}_0(r,l,-1,q)$ and a basic computation shows that,
among them, only $(q-1)(q-2)^l$ satisfy condition $(C_0)$. 
\par In this statement, note that we only consider the case
$a=-1$. The general case is derived from the following argument. Let $\alpha \in \mathcal{E}_0(r,l,a,q)$ be defined
by  
$$\alpha=[\lambda_1T,\dots,\lambda_lT,\alpha_{l+1}] \quad \text{ and}\quad \alpha^r=\epsilon_1P_a\alpha_{l+1}+\epsilon_2Q_a.$$ 
Let $v$, in $\F_q$ or $\F_{q^2}$, be such that $v^2=-a$, then define
$\beta(T)=v\alpha(vT)$. One can show that $\beta$ belongs to $\F(q)$
and that it is
defined by  
$$\beta=[-a\lambda_1T,\lambda_2T,-a\lambda_3T,\dots,a(l)\lambda_lT,\beta_{l+1}] \quad
\text{ and}\quad
\beta^r=\epsilon'_1P_{-1}\beta_{l+1}+\epsilon'_2Q_{-1},$$ 
where $\epsilon'_2=a^{r-1}\epsilon_2$, $\epsilon'_1=a^{r-1}a(l)\epsilon_1$ and $a(l)=1$ if $l$ is even or $-a$ if $l$ is odd. Consequently, we have $\beta \in \mathcal{E}_0(r,l,-1,q)$ and there is a one to one correspondance
between  the sets $\mathcal{E}_0(r,l,-1,q)$ and
$\mathcal{E}_0(r,l,a,q)$. Accordingly, condition $(C_0)$ can easily be generalized, and we have :
\vskip 0.5 cm 
\noindent {\it{If $\alpha \in \mathcal{E}_0(r,l,a,q)$ and 
  $(C_0)$ : $[-a\lambda_1,\lambda_2,\dots,a(l)(\lambda_l+\epsilon_1/\epsilon_2)]^r=a^{r-1}\epsilon_2$,
  then   $\alpha \in \mathcal{E}^*_0(r,l,a,q)$.}} 
\par The present work is organized as follows. In the next section, we introduce the basic results concerning continued fractions, used in this note. In section 3, we establish an important property of the pair $(P_a,Q_a)$ which plays a key role in the description of our particular algebraic continued fractions. In section 4, we give the proof of the theorem which is stated here below. In a short and last section, we make further comments and we give an orientation toward further studies. 
\vskip 0.5 cm
\noindent {\bf{Theorem.}}{\emph{  Let $p$ be an odd prime number, $q=p^s$, $r=p^t$, with integers $s,t\geq 1$. Let $l\geq 1$ be an integer. Let $(a,\epsilon_1,\epsilon_2)
\in (\F^*_q)^3$ be given. Let $P_a,Q_a\in \F_q[T]$ be defined as above. Let $\alpha =[a_1,a_2,\dots,a_n,\dots] \in \F(q)$ be the infinite continued
  fraction defined by
$$(a_1,\dots,a_l)=(\lambda_1T+\mu_1,\dots ,\lambda_lT+\mu_l) \quad , \text{where }\quad(\lambda_i,\mu_i)\in \F_q^*\times
\F_q\quad \text{ for }\quad 1\leq i\leq l,$$
and
$$\alpha^r=\epsilon_1P_a\alpha_{l+1}+\epsilon_2Q_a.$$
We assume that the $(2l+1)$-tuple
$(\lambda_1,\dots,\lambda_l,\mu_1,\dots,\mu_{l-1},a,\epsilon_2)$
is such that, for $1\leq i\leq l$, we can define the pair $(\delta_i,\nu_i)\in \F_q^*\times
\F_q$ in the following way :
$$\delta_1=a\lambda_1^r+\epsilon_2, \quad \nu_1=0, \quad
\text{ and }\quad \text{ for }\quad 1\leq i\leq l-1$$
$$\delta_{i+1}=a\lambda_{i+1}^r-\frac{\delta_i}{a^{r-2}\delta_i^2+(\nu_i-\mu_i^r)^2}\quad
\text{ and }\quad
\nu_{i+1}=\frac{(\nu_i-\mu_i^r)}{a^{r-2}\delta_i^2+(\nu_i-\mu_i^r)^2}.\eqno{(C_1)}$$
We also assume that the pair $(\epsilon_1,\mu_l)$ is such that we have : $$(C_2)\quad \delta_l=-a(\epsilon_1/\epsilon_2)^r\quad
\text{ and }\quad (C_3)\quad \mu_l^r=\nu_l.$$
 Then we have $\alpha \in \mathcal{E}^*(r,l,a,q)$. Moreover the sequence
of partial quotients, defined by $a_n=\lambda_nT+\mu_n$ for $n\geq 1$, is described as follows.
\newline For $n\in \N^*$, we set $f(n)=nr+l+1-r$ and $g(n)=nr+l+(1-r)/2$. We introduce the following subsets of $\N^*$: $I=\lbrace i\in \N \vert 1\leq i\leq l \rbrace$, $I^*=\lbrace i\in I\vert \nu_i-\mu_i^r \neq 0 \rbrace$, $F=\lbrace f^m(i)\quad \vert  m\geq 1 \text{ and } i\in I \rbrace$ and $G=\lbrace g^m(i)\quad \vert m\geq 1 \text{ and } i\in I^* \rbrace$. Note that the subsets $I^*$ and $G$ may both be empty: namely if $\mu_i=0$ for $1\leq i \leq l$.
\newline For $n>l+1$, we define $C(n)$ by : 
$C(n)=4a^{-1}$ if $n\notin F\cup (F+1)\cup (G+1)$ and
\begin{eqnarray*}
C(n)&=& 2a^{-1}(1-a^{-1}\lambda_i^{-r}\delta_i)^{-r^{m-1}}  \quad \text{if} \quad n=f^m(i)\quad \text{for} \quad m\geq 1\quad \text{and} \quad i\in I \quad (n\neq f(1)),\\
C(n)&=& 2a^{-1}(a\lambda_i^r\delta_i^{-1})^{r^{m-1}}  \quad \text{if} \quad n=f^m(i)+1\quad \text{for} \quad m\geq 1\quad \text{and} \quad i\in I,\\
C(n)&=& 4a^{-1}(1+a^{2-r}(\nu_i-\mu_i^r)^2\delta_i^{-2})^{-r^{m-1}}  \quad \text{if} \quad n=g^m(i)+1\quad \text{for} \quad m\geq 1\quad \text{and} \quad i\in I^*.\\
\end{eqnarray*}
Then the sequence $(\lambda_n)_{n\geq1}$ in $\F_q^*$ is defined recursively, for $n\geq l+1$, by
$$ \lambda_{l+1}=\lambda_1^r\epsilon_1^{-1}\quad \text{and} \quad \lambda_n= C(n)\lambda_{n-1}^{-1}\quad \text{for} \quad n>l+1.$$
While the sequence $(\mu_n)_{n\geq1}$ in $\F_q$ is defined, for $n\geq l+1$, as follows. 
\newline If $n\notin G\cup (G+1)$, then $\mu_n=0$. If $n=g^m(i)$ for $m\geq 1$ and $i\in I^*$, then 
$$\mu_n\lambda_n^{-1}=-\mu_{n+1}\lambda_{n+1}^{-1}=(-a)^{(r^{m-1}(2-r)+1)/2}((\nu_i-\mu_i^r)\delta_i^{-1})^{r^{m-1}}/2.$$}}

\noindent{\bf{Remark:}} If $l=1$, we observe that conditions $(C_1)$, $(C_2)$ and
$(C_3)$ reduce to $a\lambda_1^r+\epsilon_2=-a(\epsilon_1/\epsilon_2)^r$ with $\mu_1=0$. Consequently, the corresponding continued fraction
belongs to $\mathcal{E}^*_0(r,1,a,q)$. As we already observed, the example introduced at the begining of this section belongs to  $\mathcal{E}^*_0(p,1,-1,p)$ and we have the desired condition $\delta_1=-1+\epsilon_2=\epsilon_1/\epsilon_2$. Besides, the reader may check that the description of the continued fraction given in this example can be derived from the formulas stated in this theorem.
\newline Also, in the simplest case
$(\mu_1,\mu_2,\dots,\mu_l)=(0,0,\dots,0)$, i.e. $\alpha \in \mathcal{E}_0(r,l,a,q)$, $(C_1)$ implies  inductively $\nu_i=0$ for $1\leq i\leq l$. Consequently, $(C_1)$ reduces to
$\delta_{i+1}=a\lambda_{i+1}^r-a^{2-r}\delta_i^{-1}$ for $1\leq i\leq
l-1$. One can check that conditions $(C_1)$, $(C_2)$ and
$(C_3)$ then reduce to the sufficient condition $(C_0)$, already stated above, in order to have $\alpha \in
\mathcal{E}^*_0(r,l,a,q)$. 
\par Before concluding this section, we make a more general comment on the subject discussed in this note. This comment underlines the place taken by the family of power series described here among the hyperquadratic power series. All the known examples of algebraic continued fractions, in odd characteristic, having all partial quotients of degree 1, are related to continued fractions generated in the way presented above. To be more precise, let us denote by  $\mathcal {E^*}(q)$ the union of the sets $\mathcal{E^*}(r,l,a,q)$, for all $l\geq 1$, all $a\in \F_q^*$ and all $r=p^t$, with $t\geq 0$. Then we conjecture that, if $\alpha \in \H(q)$ ($q$ odd) and if all its partial quotients are of degree 1, then there is a linear fractional transformation $f(x)=(ax+b)/(cx+d)$, with $(a,b,c,d) \in \F_q[T]^4$ and $ad-bc\in \F_q^*$, and $\beta \in \mathcal {E^*}(q)$, such that $\alpha(T)=f(\beta(\lambda T+\mu))$ where $(\lambda,\mu)\in \F_q^*\times \F_q$.
\newline Also, it is a known fact that if $\alpha \in \mathcal {E^*}(q)$ and $f(x)=(ax+b)/(cx+d)$, with $(a,b,c,d) \in \F_q[T]^4$ is a linear fractional transformation, then $f(\alpha)$ belongs to $\H(q)$. Moreover $\alpha$ and $f(\alpha)$ have the same algebraic degree, and $f(\alpha)$ has also bounded partial quotients. Hence the set $\mathcal {E^*}(q)$ generates a subset of elements in $\H(q)$ with bounded partial quotients. However, most elements in $\H(q)$ have unbounded partial quotients. The reader may see the introduction of $[8]$, for more information on this matter.
\newline Finally, thinking of a famous conjecture in number theory in the classical context of real numbers, we ask the following question : are there algebraic irrational power series, in odd characteristic, which are not hyperquadratic and which have partial quotients of bounded degrees ?

\vskip 0.5 cm
\noindent{\bf{2. Notation and basic formulas for continued fractions}}
\par Let $W=w_1,w_2,\ldots,w_n$ be a sequence of variables over a ring $\mathbb A$. We set $\vert W\vert =n$ for the length of the word $W$.
We define the following operators for the word $W$.
\begin{align*}
W'&=w_2, w_3,\ldots, w_n \quad \text{ or }\quad  W'=\emptyset \quad
\text{ if } \quad \vert W\vert =1.\\
W''&= w_1, w_2, \ldots, w_{n-1}\quad  \text{ or} \quad W''=\emptyset
\quad \text{if } \quad \vert W\vert =1.\\
W^*&=w_n, w_{n-1},\ldots, w_1.
\end{align*}
We consider the finite continued fraction associated to $W$ to be
$$[W]=[w_1,w_2,\ldots,w_n]=w_1+\cfrac{1}{w_2+\cfrac{1}{\ddots+\cfrac{1}{w_n}}}.$$
This continued fraction is a quotient of multivariate polynomials,
usually called continuants, built upon the variables $w_1,w_2,\dots,w_n$. 
More details about these polynomials can
be found, for example in [13], and also in [6] (although here,
trying to simplify, we adopt different notation).
The continuant buildt on $W$ will be denoted $<W>$. We now recall the
definition of this sequence of multivariate polynomials. 
\newline Set $<\emptyset>=1$. If the sequence $W$ has only one element, then we have $<W>=W$. 
Hence, with the above notation, the continuants can be computed,
recursively on the lenght $\vert W\vert$, by the following
formula 
$$<W>=w_1<W'>+<(W')'> \quad \text{ for} \quad \vert W \vert \geq 2. \eqno{(1)}$$
Thus, with this notation, for any finite word $W$, the finite continued fraction $[W]$ satisfies
$$[W]=\frac{<W>}{<W'>}.$$
It is easy to check that the polynomial $<W>$ is, in a certain sense, symetric in the variables $w_1,w_2,\ldots,w_n$. 
Hence we have $<W^*>=<W>$ and this symetry implies the classical formula 
$$[W^*]=\frac{<W>}{<W''>}.$$
The continuants satisfy a number of useful identities. First we will need a generalization of (1). 
For any finite sequences $A$ and $B$, of variables over $\mathbb A$, defining $A,B$ 
as the concatenation of sequences $A$ and $B$, we have
$$<A,B>=<A><B>+<A''><B'>.\eqno{(2)}$$
Secondly, using induction on $\vert W \vert$, we have the
following classical identity 
$$<W><(W')''>-<W'><W''>=(-1)^{ \vert W \vert}\quad \text{ for} \quad \vert W \vert \geq 2.\eqno{(3)}$$
 Now, let $y$ be an invertible element of $\mathbb A$, then we define $y
\cdot W$ as the following sequence
$$y\cdot W = y w_1, y^{-1}w_2,\ldots, y^{(-1)^{n-1}}w_n.$$
With these notations, it is easy to check that we have $y[W]=[y \cdot
W]$ and more precisely
$$<y \cdot W>=<W> \quad \text{ if} \quad \vert W \vert \text{ is
  even}\quad \text{ and } \quad <y \cdot W>=y<W>\text{ if} \quad \vert
W \vert \text{ is odd} . \eqno{(4)}$$
 Let us come back to the notation used in the introduction. If
$\alpha \in \F(q)$ is irrational (rational), then it can be expanded as an infinite (finite) continued fraction $\alpha=[a_1,a_2,\dots,a_n,\dots]$, where the $a_n\in \F_q[T]$ are called the partial quotients. We have $\deg(a_n)>0$ for $n>1$. For $n\geq 1$, we
set $x_n=<a_1,a_2,\dots,a_n>$ and $y_n=<a_2,\dots,a_n>$, with $x_0=1$ and $y_0=0$. The rational $x_n/y_n=[a_1,a_2,\dots,a_n]$ is called a
convergent to $\alpha$. The continued fraction expansion of an irrational element measures the quality of its rational approximation. The convergents of $\alpha$ are the best rational approximations and we have $\vert\alpha-x_n/y_n\vert=\vert a_{n+1}\vert^{-1}\vert y_n\vert^{-2}$. If the partial quotients have bounded degrees then the element is said to be badly approximable. Let us recall that we also have
$\alpha=(x_n\alpha_{n+1}+x_{n-1})/(y_n\alpha_{n+1}+y_{n-1})$ for
$n\geq 1$, where $\alpha_{n+1}=[a_{n+1},a_{n+2},\dots]$ is the
complete quotient. 
\par We will also make use of the following general and basic lemma. 
\newline {\bf{Lemma 0. }}{\emph{ Let $W$ be a finite word, with $\vert W \vert\geq 2$. Let $a$ be a variable over $\mathbb A$, then we have $$[W]+a=[W,b],$$ where
$b=(-1)^{ \vert W \vert-1}<W'>^{-2}a^{-1}-<(W')''><W'>^{-1}.$}}
\newline Indeed, we have
$$[w_1,w_2,\dots,w_n,b]=\frac{x_nb+x_{n-1}}{y_nb+y_{n-1}}=\frac{x_n}{y_n}+\frac{(-1)^{n-1}}{y_n(y_nb+y_{n-1})}=[w_1,\dots,w_n]+a.$$
A first publication of the idea under this statement is due to M. Mend\`es France (see [11, p. 209]).
\vskip 0.5 cm
\noindent{\bf{3. A finite continued fraction in $\F_q(T)$}}
\par We recall that $p$ is an odd prime number, $q=p^s$ and $r=p^t$ where $s$ and $t$ are positive integers. For $a\in \F_q^*$, we consider the polynomials $P_a$ and $Q_a$ in $\F_q[T]$  defined by :
$$P_a(T)=(T^2+a)^{(r-1)/2} \quad \text{ and } \quad Q_a(T)=a^{-1}(TP_a(T)-T^r).\eqno{(5)}$$
The following sequence $(F_n)_{n\geq 0}$ of polynomials in $\F_p[T]$ was introduced by Mills and Robbins [12, p. 400] (see also [6, p. 331]). This sequence is defined recursively by
$$F_0=1, \quad F_1=T \quad \text{ and } \quad F_n=TF_{n-1}+F_{n-2}\quad \text{ for } \quad n\geq 2.\eqno{(6)}$$
From $(6)$, we clearly have the finite continued fraction expansion
$$F_n/F_{n-1}=[T,T,\dots,T] \quad (n \quad \text{terms}).$$
This sequence can be regarded as the analogue in the function field case of the Fibonacci sequence of integers. By elementary computations (see [6, pp. 331-332]), one can check that the following formulas hold in $\F_p[T]$:
$$F_{r-1}=P_4  \quad \text{ and } \quad  F_{r-2}=-2Q_4.$$
Consequently, with $(4)$, we can write
$$P_4/Q_4=-2F_{r-1}/F_{r-2}=[-2T,-T/2,\dots,-2T,-T/2] .\eqno{(7)}$$
Now, let $v\in \F_{q^2}$ be such that $v^2=a/4$. From (5) we get  
$$P_a(T)=(a/4)^{(r-1)/2}P_4(T/v) \quad \text{ and } \quad  Q_a(T)=v(a/4)^{(r-3)/2}Q_4(T/v).\eqno{(8)}$$
Therefore, by $(7)$ and $(8)$, we have 
$$P_a/Q_a=(a/4v)(P_4/Q_4)(T/v)=v[-2T/v,-T/2v,\dots,-2T/v,-T/2v] $$
and, with $(4)$, finally
$$P_a/Q_a=[-2T,-2T/a,\dots,-2T,-2T/a] \quad (r-1 \quad \text{terms}).\eqno{(9)}$$
\par Let us make a remark on the infinite continued fraction $\omega=[T,T,\dots,T,\dots]$ in $\F(p)$. This element is quadratic and it
clearly satisfies $\omega^2=T\omega+1$ (it is an analogue of the golden
mean in the case of real numbers). One can prove, for
all $n\geq 1$, the equality $\omega^{n+1}=F_n\omega+F_{n-1}$. Consequently, we obtain $\omega^r=P_4\omega-2Q_4$. Since we have  
$\omega=\omega_{l+1}$, it follows that $\omega$ belongs to
$\mathcal{E}^*_0(r,l,4,p)$, for all $r$, all $p>2$ and all $l\geq
1$. It is well known that, also in the case of power series over a finite field, quadratic continued
fractions are characterized by an ultimately periodic sequence of
partial quotients. For a general element in  $\mathcal{E^*}(r,l,a,q)$,
this sequence is not so and therefore this element is not quadratic. However
the precise algebraic degree of such an element is generally unknown. Concerning this matter, we can make an observation about the example introduced at the begining of this article. Indeed, one could check that, for the particular value $\epsilon=1/2$, this element $\alpha$ is actually quadratic and we have, for all $p\geq 3$, $\alpha(T)=(\sqrt{-1}/2)\omega ((2/\sqrt{-1})T)$. While if $\epsilon\neq 1/2$ the algebraic degree of the corresponding element might well be $p+1$. 
\par The aim of the following proposition is to give a generalization for the continued fraction expansion $(9)$ concerning the pair $(P_a,Q_a)$.
\vskip 0.5 cm 
\noindent {\bf{Proposition 1. }}{\emph{ Let $r$, $q$, $a$, $P_a$ and $Q_a$ be defined as above. We set 
$k=(r-1)/2$. Let $x\in \F_q$, we set $\omega=1+a^{2-r}x^2$. Then $P_a$ and $Q_a+x$ are coprime polynomials in $\F_q[T]$ if and only 
if $\omega \neq 0$. We assume that $\omega \neq 0$. We define $2k$ polynomials in $\F_q[T]$ :
$$ v_i=-2T\quad \text{ and } \quad v_{k+1+i}=-\omega^{(-1)^{i+1}}2T \quad \text{
  for } \quad 1\leq i \leq k-1,$$
$$v_k=-2T-(-a)^{1-k}x\quad  \text{ and } \quad v_{k+1}(x)=\omega^{-1}(-2T+(-a)^{1-k}x).$$
We set $W(a,x)=v_1,v_2/a,v_3,v_4/a,\ldots,v_{2k-1},v_{2k}/a$. Then we have the following equalities in $\F_q(T)$:
 $$P_a(Q_a+x)^{-1}=[W(a,x)],\quad <W(a,x)>=\omega(k)a^{-k}P_a ,\quad <W'(a,x)>=\omega(k)a^{-k}(Q_a +x)$$ 
 and also 
$$<W''(a,x)>=\omega(k)\omega^{(-1)^{k+1}}a^{1-k}(Q_a -x),$$
where $\omega(k)=1$ if $k$ is even and $\omega(k)=\omega^{-1}$ if $k$ is odd.}}
\vskip 0.5 cm 
\noindent{\bf{Proof:}} We denote by $v\in \F_{q^2}$ a square
  root of $-a$, so that $\pm v$ are the only roots of $P_a$. Hence we
  see that $P_a$ and $Q_a+x$ are coprime if and only if
  $(Q_a(v)+x)(Q_a(-v)+x)\neq 0$. From $(5)$, we obtain $Q_a(\pm v)=\mp a^{-1}v^r$. Therefore, this becomes $x^2-a^{-2}v^{2r}\neq 0$. Observing that $v^{2r}=-a^r$, we obtain the desired condition. 
\newline First, since $\omega=1$ if $x=0$, we observe that $W(a,0)=-2T,-2T/a,-2T,\ldots,-2T/a$. Hence, with our notations, equality (9) can be written as
$$\frac{P_{a}}{Q_{a}}=\frac{<W(a,0)>}{<W'(a,0)>}. $$
Since the numerators of both fractions above have the same degree, $r-1$ in $T$, it follows that 
there exists $u\in \F_q^*$ such that 
$$<W(a,0)>=uP_a\quad \text{ and } \quad <W'(a,0)>=uQ_a.\eqno{(10)}$$
For a continuant built from polynomials in $T$, the leading coefficient is obtained as the leading coefficient of the product 
of its terms. Consequently, the leading coefficient of $<W(a,0)>$ is $(-2)^{r-1}a^{-k}=a^{-k}$ while $P_{a}$ is unitary. 
This implies $u=a^{-k}$. Now we observe that $k$ can be even or odd. If $p=4m+1$ then $k$ is allways even, while if $p=4m+3$ then $k$ has the same parity as $t$ if $r=p^t$. 
\newline We set $W(a,0)=a_1,a_2,\ldots,a_{r-1}$ and $W(a,x)=b_1,b_2,\ldots,b_{r-1}$. We will use the notation $a(i)=1$ if $i$ is odd and $a(i)=a^{-1}$ if $i$ is even. Hence we have $a_i=-2a(i)T$ and $b_i=a(i)v_i$ for $1\leq i\leq 2k$. To shorten the writting, we denote $<a_i,\dots,a_j>$ by $A_{i,j}$ and similarly  $<b_i,\dots,b_j>$ by $B_{i,j}$. According to these notations we have $<W(a,0)>=A_{1,2k}$ and $<W(a,x)>=B_{1,2k}$. From the definition of both sequences
$W(a,x)$ and $W(a,0)$, we have $b_i=a_i$  and  $b_{k+1+i}=\omega^{(-1)^{i+1}}a_{k+1+i}$  for $ 1\leq i \leq k-1$. 
Consequently, we get 
$$B_{1,k-1}=A_{1,k-1} \quad \text{ and } \quad B_{2,k-1}=A_{2,k-1}.\eqno{(11)}$$
But also, by $(4)$, according to the parity of $k$,  
$$B_{k+2,2k}=\omega\omega(k)A_{k+2,2k}\quad \text{ and } \quad B_{k+3,2k}=\omega(k)A_{k+3,2k} .\eqno{(12)}$$
Applying $(2)$, by (11) and (12), we can
write  
$$B_{1,2k}=B_{1,k}B_{k+1,2k}+B_{1,k-1}B_{k+2,2k}=B_{1,k}B_{k+1,2k}+\omega \omega(k)A_{1,k-1}A_{k+2,2k}\eqno{(13)}$$
$$B_{2,2k}=B_{2,k}B_{k+1,2k}+B_{2,k-1}B_{k+2,2k}=B_{2,k}B_{k+1,2k}+\omega \omega(k)A_{2,k-1}A_{k+2,2k}.\eqno{(14)}$$
For notational convenience define $x_k=(-a)^{1-k}x$. From the definition of $v_k$ and $v_{k+1}$, we have
$$b_k=a_k-a(k)x_k\quad \text{ and } \quad b_{k+1}=\omega^{-1}(a_{k+1}+a(k+1)x_k). \eqno{(15)}$$
By (2), we have $B_{1,k}=B_{1,k-1}b_k+B_{1,k-2}$. Again by (2), (11) and (15), this becomes 
$$B_{1,k}=A_{1,k-1}b_k+A_{1,k-2}=A_{1,k}-a(k)x_kA_{1,k-1}. \eqno{(16)}$$
In the same way, by (2), (11) and (15), we get
$$B_{2,k}=A_{2,k-1}b_k+A_{2,k-2}=A_{2,k}-a(k)x_kA_{2,k-1}. \eqno{(17)}$$
 By (2), (12) and (15), since $B_{k+1,2k}=B_{k+2,2k}b_{k+1}+B_{k+3,2k}$, we also get
$$B_{k+1,2k}=\omega(k)A_{k+2,2k}(a_{k+1}+a(k+1)x_k)+\omega(k)A_{k+3,2k}$$
and this becomes $$B_{k+1,2k}=\omega(k)(A_{k+1,2k}+a(k+1)x_kA_{k+2,2k}). \eqno{(18)}$$
By $(4)$, according to the parity of $k$, we also obtain 
$$A_{k+1,2k}=a(k+1)A_{1,k} \quad \text{ and }\quad A_{k+2,2k}=a(k)A_{1,k-1}.\eqno{(19)}$$
From $(19)$, we have $a(k)A_{1,k-1}A_{k+1,2k}-a(k+1)A_{1,k}A_{k+2,2k}=0$. We also have $a(k)a(k+1)x_k^2=a^{-1}a^{2-2k}x^2=a^{2-r}x^2$.
Consequently, by multiplication, from (16) and (18), we get
$$B_{1,k}B_{k+1,2k}=\omega(k)(A_{1,k}A_{k+1,2k}-a^{2-r}x^2A_{1,k-1}A_{k+2,2k}). \eqno{(20)}$$
In the same way, by multiplication, from (17) and (18), we get 
$$B_{2,k}B_{k+1,2k}=\omega(k)(A_{2,k}A_{k+1,2k}-a^{2-r}x^2A_{2,k-1}A_{k+2,2k}+X), \eqno{(21)}$$
where, according to $(19)$, using $(3)$ and $a(k)a(k+1)=a^{-1}$,  we have
$$X=a(k)a(k+1)x_k(A_{1,k-1}A_{2,k}-A_{2,k-1}A_{1,k})=a^{-1}x_k(-1)^{k-1}=a^{-k}x.\eqno{(22)}$$ 
Combining (13) and (20), using $(2)$ and since $\omega=1+a^{2-r}x^2$, we get 
$$B_{1,2k}=\omega(k)(A_{1,k}A_{k+1,2k}+(\omega-a^{2-r}x^2)A_{1,k-1}A_{k+2,2k})=\omega(k)A_{1,2k}.$$
By $(10)$, recalling that $u=a^{-k}$, this becomes  
$$<W(a,x)>=B_{1,2k}=\omega(k)A_{1,2k}=\omega(k)<W(a,0)>=\omega(k)a^{-k}P_a.$$
In the same way, combining (14), (21) and (22), we obtain
$$B_{2,2k}=\omega(k)(A_{2,k}A_{k+1,2k}+(\omega-a^{2-r}x^2)A_{2,k-1}A_{k+1,2k}+a^{-k}x)=\omega(k)(A_{2,2k}+a^{-k}x).$$
By $(10)$, with $u=a^{-k}$, this becomes  
$$<W'(a,x)>=B_{2,2k}=\omega(k)(A_{2,2k}+a^{-k}x)=\omega(k)(<W'(a,0)>+a^{-k}x)=\omega(k)a^{-k}(Q_a+x).$$ 
Consequently, we get
$$[W(a,x)]=\frac{<W(a,x)>}{<W'(a,x)>}=\frac{\omega(k)a^{-k}P_a}{\omega(k)a^{-k}(Q_a+x)}=P_a(Q_a+x)^{-1}.$$
Moreover, from the definition of the sequence $W(a,x)$, we 
observe the "pseudo-symetry" between $W(a,x)$ and $W(a,-x)$, i.e. $W(a,x)=a\omega^{(-1)^{k+1}}\cdot W^*(a,-x)$. Finally, 
using this equality, by $(4)$ and since $\vert
W'(a,-x) \vert=2k-1$ is odd,  we obtain
$$<W''(a,x)>=<W''^*(a,x)>=<a\omega^{(-1)^{k+1}}\cdot W'(a,-x)>=\omega(k)\omega^{(-1)^{k+1}}a^{1-k}(Q_a-x).$$ 
So the proof of Proposition 1 is complete.
\vskip 0.5 cm
\noindent{\bf{4. Proof of the theorem}}

\par Throughout this section, the integers $p$, $q$ and $r$, as well as $a\in \F_q^*$ and $P_a$, $Q_a$ in $\F_q[T]$, are defined as above. Moreover, as above, we set $k=(r-1)/2$. We need the following lemma, which is a straightforward consequence of Lemma 0 from Section 2 and of Proposition 1 from Section 3. 
\newline {\bf{Lemma 1.}}{\emph{ Let $b_0\in \F_q[T]$ and $y\in \F_q^*$. For
  $a\in \F_q^*$ and $x\in \F_q$, assuming that $\omega=1+a^{2-r}x^2\neq 0$, as above we denote by $W(a,x)$ the
  sequence of the $r-1$ partial quotients of the rational function $P_a(Q_a+x)^{-1}$. 
Then, for $X\in \F(q)$, we have the formal identity:
$$[b_0,y\cdot W(a,x)]+X=[b_0,y\cdot W(a,x),Y],$$
where 
$$Y=\omega^{(-1)^{k+1}} (\omega a^{r-1}P_a^{-2}X^{-1}-ya(Q_a-x)P_a^{-1}).$$}}
\noindent{\bf{Proof:}} According to  Lemma 0 in Section 2, we can write
$$[b_0,y\cdot W(a,x)]+X=[b_0,y\cdot W(a,x),Y],$$
where $Y$ is linked to $X$ as follows
$$Y=(-1)^{r-1}<y\cdot W(a,x)>^{-2}X^{-1}-<(y\cdot W(a,x))''><y\cdot W(a,x)>^{-1}.\eqno{(23)}$$
We recall that 
$<y\cdot W(a,x)>=<W(a,x)>$, since the sequence of terms is of even length $r-1$. In the same way,
since $r-2$ is odd, we also have  $<(y\cdot W(a,x))''>=y<W''(a,x)>$.
Applying Proposition 1, we have $<W(a,x)>=\omega(k)a^{-k} P_a$ and
$<W''(a,x)>=\omega(k)\omega^{(-1)^{k+1}}a^{1-k}(Q_a-x)$. Consequently,
(23) becomes
$$Y=\omega(k)^{-2}a^{2k}P_a^{-2}X^{-1}-ya \omega^{(-1)^{k+1}} (Q_a-x)P_a^{-1}.\eqno{(24)}$$
Since $\omega(k)^{-2}=\omega^{1+(-1)^{k+1}}$ and $2k=r-1$, (24) implies the conclusion of this lemma.
\par The proof of the theorem relies on the following proposition. 
\newline {\bf{Proposition 2.}}{\emph{ Let $p$, $q$ and $r$ be as above. Let $\alpha=[a_1,a_2,\ldots,a_n,\ldots]$ be an irrational element of
$\F(q)$. For an integer $n\geq 1$, we set $f(n)=(n-1)r+l+1$. For an index $n\geq 1$, we assume that $a_n=\lambda_nT+\mu_n$, where $(\lambda_n,\mu_n)\in \F_q^*\times \F_q$
and that $\alpha_n$ and $\alpha_{f(n)}$ are linked by the following equality
$$\alpha_n^r=\epsilon_{1,n} P_a \alpha_{f(n)}+\epsilon_{2,n}Q_a+\nu_n,$$
where $(\epsilon_{1,n},\epsilon_{2,n},\nu_n) \in (\F_q^*)^2\times
\F_q$. We set $\delta_n=a\lambda_n^r+\epsilon_{2,n}$. \newline First we assume
that $\delta_n\neq 0$. We set $\pi_n=(\nu_n-\mu_n^r)\delta_n^{-1}$ and
$\omega_n=1+a^{2-r}\pi_n^2$. We assume that $\omega_n\neq 0$. The word $W(a,\pi_n)$ is defined in Proposition 1. Then we have
 $$a_{f(n)}=\epsilon_{1,n}^{-1}\lambda_n^rT \quad \text{ and } \quad
a_{f(n)+1},\ldots,a_{f(n)+r-1}=(-\epsilon_{1,n} \delta_n^{-1})\cdot W(a,\pi_n).$$
Moreover we have
$\alpha_{n+1}^r=\epsilon_{1,n+1}P_a\alpha_{f(n+1)}+\epsilon_{2,n+1}Q_a+\nu_{n+1}$, where 
\begin{eqnarray*}
\epsilon_{1,n+1}&=& a^{1-r}\epsilon_{1,n}^{-1} \omega_n^{(-1)^k-1}\\
\epsilon_{2,n+1}&=& -a^{2-r}(\omega_n\delta_n)^{-1}\\
\nu_{n+1}&=& a^{2-r}(\nu_n-\mu_n^r)\omega_n^{-1}\delta_n^{-2}.\\
\end{eqnarray*}
Finally, if $\delta_n=0$ then we have $a_{f(n)}=\epsilon_{1,n}^{-1}\lambda_n^rT$, but $\deg(a_{f(n)+1})>1$.
}}
\vskip 0.5 cm
\noindent{\bf{Proof of Proposition 2:}} By hypothesis , we have $\alpha_n=[\lambda_nT+\mu_n, \alpha_{n+1}]$ and  also
$$\alpha_n^r=\epsilon_{1,n} P_a \alpha_{f(n)}+\epsilon_{2,n}Q_a+\nu_n.\eqno{(25)}$$
Therefore, combining the first equality with (25), and since $\alpha_n^r=[a_n^r,\alpha_{n+1}^r]$, we can write 
$$[\lambda_n^rT^r+\mu_n^r-\epsilon_{2,n}Q_a-\nu_n,\alpha_{n+1}^r]=\epsilon_{1,n}
P_a\alpha_{f(n)}.\eqno{(26)}$$
Recalling that $P_a$ and $Q_a$ satisfy the equality $T^r=TP_a-aQ_a$,
we obtain, with our notation,  
$$\lambda_n^rT^r-\epsilon_{2,n}Q_a=\lambda_n^rTP_a-\delta_nQ_a.\eqno{(27)}$$
Combining (26) and (27), we get 
$$[\frac{\lambda_n^rTP_a-\delta_nQ_a+\mu_n^r-\nu_n}{\epsilon_{1,n} P_a},
\epsilon_{1,n} P_a\alpha_{n+1}^r]=\alpha_{f(n)}.\eqno{(28)}$$
Assuming that $\delta_n \neq 0$, with our notation, since $\pi_n=(\nu_n-\mu_n^r)\delta_n^{-1}$, $(28)$ can be written as
$$\epsilon_{1,n}^{-1}\lambda_n^rT-\epsilon_{1,n}^{-1}\delta_nP_a^{-1}(Q_a+\pi_n)+
\epsilon_{1,n}^{-1}P_a^{-1} \alpha_{n+1}^{-r}=\alpha_{f(n)}.\eqno{(29)}$$
Applying Proposition 1, we have $P_a(Q_a+\pi_n)^{-1}=[W(a,\pi_n)]$. We set $y=-\epsilon_{1,n}\delta_n^{-1}$. Then we have
$$-\epsilon_{1,n}\delta_n^{-1}P_a(Q_a+\pi_n)^{-1}=y[W(a,\pi_n)]=[y\cdot W(a,\pi_n)].$$
 Consequently, if we set $b_0=\epsilon_{1,n}^{-1}\lambda_n^rT$ and $X=\epsilon_{1,n}^{-1}P_a^{-1} \alpha_{n+1}^{-r}$, (29) becomes 
$$b_0+\frac{1}{[y\cdot W(a,\pi_n)]}+X=\alpha_{f(n)},\eqno{(30)}$$
which is $[b_0,y\cdot W(a,\pi_n)]+X=\alpha_{f(n)}$. Since $\omega_n=1+a^{2-r}\pi_n^2\neq 0$, we can apply Lemma 1 above. We get
$$[\epsilon_{1,n}^{-1}\lambda_n^rT,(-\epsilon_{1,n}\delta_n^{-1})\cdot W(a,\pi_n),Y]=\alpha_{f(n)}.\eqno{(31)}$$
This lemma gives
$$Y=\omega_n^{(-1)^{k+1}} \epsilon_{1,n}P_a^{-1}(\omega_na^{r-1}\alpha_{n+1}^r+a\delta_n^{-1}(Q_a-\pi_n)).\eqno{(32)}$$
We have $\vert Y\vert =\vert P_a^{-1}\alpha_{n+1}^r\vert >1$, consequently (31) implies 
$$a_{f(n)}=\epsilon_{1,n}^{-1}\lambda_n^rT\quad \text{ and}\quad  a_{f(n)+1},\dots,
a_{f(n)+r-1}=(-\epsilon_{1,n}\delta_n^{-1})\cdot W(a,\pi_n).$$
 But also
$Y=\alpha_{f(n)+r}=\alpha_{f(n+1)}$. Hence, by (32),
we get
$$\omega_n^{(-1)^{k}}\epsilon_{1,n}^{-1}P_a\alpha_{f(n+1)}= \omega_na^{r-1}\alpha_{n+1}^r+a\delta_n^{-1}(Q_a-\pi_n)\eqno{(33)}$$
and (33) becomes
$$\alpha_{n+1}^r=a^{1-r}\epsilon_{1,n}^{-1}\omega_n^{(-1)^{k}-1}P_a\alpha_{f(n+1)}-a^{2-r}\delta_n^{-1}\omega_n^{-1}
Q_a+a^{2-r}\pi_n\delta_n^{-1}\omega_n^{-1}.\eqno{(34)}$$
From (34), we obtain the desired formulas for $\epsilon_{1,n+1}$,
$\epsilon_{2,n+1}$ and $\nu_{n+1}$. 
\newline Finally if $\delta_n=0$, from $(28)$, we get
$$\alpha_{f(n)}=\epsilon_{1,n}^{-1}\lambda_n^rT+\frac{(\mu_n^r-\nu_n)\alpha_{n+1}^r+1}{\epsilon_{1,n} P_a\alpha_{n+1}^r}=\epsilon_{1,n}^{-1}\lambda_n^rT+Z.$$
We have $\vert Z\vert <\vert T\vert^{-1}$, consequently we get $a_{f(n)}=\epsilon_{1,n}^{-1}\lambda_n^rT$ and $\deg(a_{f(n)+1})=\deg(Z^{-1})>1$.
\newline So the proof of Proposition 2 is complete.
\vskip 0.5 cm

\noindent{\bf{Proof of the theorem :}} We start from $\alpha \in \mathcal{E}(r,l,a,q)$, satisfying $$\alpha^r=\epsilon_1P_a\alpha_{l+1}+\epsilon_2Q_a. \eqno{(I_1)}$$
Recalling that $a_1=\lambda_1 T+\mu_1$, we set
 $$\delta_1=a\lambda_1^r+\epsilon_2 \quad (DL_1)\quad \text{ and  } \quad \nu_1=0.  \quad (N_1)$$
In the sequel from the triple $(\delta_n,\nu_n,\mu_n)$ in $\F_q^3$, if $\delta_n \neq 0$, as above, we define $\pi_n=(\nu_n-\mu_n^r)\delta_n^{-1}$ and
$\omega_n=1+a^{2-r}\pi_n^2$. By $(C_1)$,  or by $(C_2)$ and $(C_3)$ if $l=1$, we have $\delta_1 \neq 0$ and $a^{r-2}\delta_1^2+(\nu_1-\mu_1^r)^2 \neq 0$. Therefore we have $\delta_1 \omega_1\neq 0$ and we can apply Proposition 2. Hence, with $f(1)=l+1$, we get $r$ partial quotients, from $a_{l+1}$ to $a_{l+r}$, all of degree 1. The following equality holds
$$\alpha_2^r=\epsilon_{1,2}P_a\alpha_{f(2)}+\epsilon_{2,2}Q_a+\nu_2,\eqno{(I_2)}$$
 where $\epsilon_{1,2}$, $\epsilon_{2,2}$ and $\nu_2$ are as stated in Proposition 2. Observe that $a_2=\lambda_2 T+\mu_2$ if $l>1$, but also if $l=1$. Indeed, if $l=1$, then $a_2=a_{f(1)}=\epsilon_1^{-1}\lambda_1^rT$. Consequently, we can consider $\delta_2$ and we have 
$$\delta_2=a\lambda_2^r-a^{2-r}(\omega_1\delta_1)^{-1} \quad (DL_2)\quad \text{ and  } \quad \nu_2=a^{2-r}(\nu_1-\mu_1^r)\omega_1^{-1}\delta_1^{-2}.  \quad (N_2)$$ 
 If $l=1$, we might have $\delta_2=0$ and this would imply $\alpha \notin \mathcal{E}^*(r,l,a,q)$. However, if $l=1$ we will see here below that $\delta_2=a^{1-r}(\delta_1\epsilon_1^{-1})^r\neq 0$. If $l>1$, again by $(C_1)$, or by $(C_2)$ and $(C_3)$ if $l=2$, we have $\delta_2\neq 0$ and $\omega_2=1+a^{2-r}\pi_2^2\neq 0$. Consequently, Proposition 2 can be applied again. This process can be carried on as long as  we have $a_n=\lambda_n T+\mu_n$ and $\delta_n\omega_n\neq 0$. As long as this process carries on, $\delta_{n}$ and $\omega_n$ are defined by means of the following recursive formulas :
 $$\delta_{n}=a\lambda_{n}^r-a^{2-r}(\omega_{n-1}\delta_{n-1})^{-1} \eqno{(DL_n)}$$
and for $\omega_n$, since $\omega_n=1+a^{2-r}(\nu_n-\mu_n^r)^2\delta_n^{-2}$, via
$$\nu_{n}=a^{2-r}(\nu_{n-1}-\mu_{n-1}^r)\omega_{n-1}^{-1}\delta_{n-1}^{-2}=a^{2-r}\pi_{n-1}(\omega_{n-1}\delta_{n-1})^{-1}.\eqno{(N_n)}$$ 
At each stage, we have $$\alpha_n^r=\epsilon_{1,n} P_a \alpha_{f(n)}+\epsilon_{2,n}Q_a+\nu_n,\eqno{(I_n)}$$ where $\epsilon_{2,n}=\delta_{n}-a\lambda_{n}^r$. While $\epsilon_{1,n}$ is defined recursively by $\epsilon_{1,1}=\epsilon_1$ and
 $$\epsilon_{1,n}=a^{1-r}\epsilon_{1,n-1}^{-1} \omega_{n-1}^{(-1)^k-1} .\eqno{(E_n)}$$
Moreover, also by Proposition 2, we have $a_m=\lambda_mT+\mu_m$, for $f(n)\leq m < f(n+1)$. Let us describe these partial quotients. We recall the notation $a(i)=1$ if $i$ is odd and  $a(i)=a^{-1}$ if $i$ is even. Combining Proposition 1 and Proposition 2, the following equalities hold : 
$$\lambda_{f(n)}=\epsilon_{1,n}^{-1}\lambda_n^r \quad (L_{n,0}), \quad \lambda_{f(n)+i}=2a(i)(\delta_n\epsilon_{1,n}^{-1})^{(-1)^i}\quad \text{for} \quad 1\leq i\leq k \eqno{(L_{n,i})}$$
$$ \text{and} \quad \lambda_{f(n)+i}=2a(i)(\delta_n\epsilon_{1,n}^{-1})^{(-1)^i}\omega_n^{(-1)^{k-i}}\quad \text{for} \quad k+1\leq i\leq r-1.\eqno{(L_{n,i})}$$
And also 
$$\mu_{f(n)+i}=0\quad \text{ for} \quad 0\leq i\leq r-1\quad \text{ and} \quad i\neq k,k+1.\eqno{(M_{n,i})}$$
$$\mu_{f(n)+k}=(-a)^{1-k}\lambda_{f(n)+k}\pi_n/2.\eqno{(M_{n,k})}$$
$$\mu_{f(n)+k+1}=-\mu_{f(n)+k}\lambda_{f(n)+k+1}\lambda_{f(n)+k}^{-1}.\eqno{(M_{n,k+1})}$$
From the equalities $(L_{n,i})$, by multiplication, we easily get the following equalities for $1\leq i\leq r-2$
$$\lambda_{f(n)+k}\lambda_{f(n)+k+1}=4a^{-1}\omega_n^{-1}\quad (X_{n,k})\quad \text{ and for}\quad i\neq k \quad   \lambda_{f(n)+i}\lambda_{f(n)+i+1}=4a^{-1}.\eqno{(X_{n,i})}$$
Our aim is to show that the quantities $\delta_n$ and $\omega_n$ can be defined, through the recursive formulas $(DL_n)$ and $(N_n)$, up to infinity. That is to say that we have $\delta_n\omega_n\neq 0$ at each stage. The first hypothesis of the theorem, namely $(C_1)$, imply that we can define recursively, by the above formulas $(DL_i)$ and $(N_i)$, $\delta_i$ and $\omega_i$ in $\F_q^*$, for $i=1,\dots,l-1$. Conditions $(C_2)$ and $(C_3)$, will turn out to be important to keep the process going on. Now, by $(C_2)$ and $(C_3)$, we observe that we also have $\delta_l\neq 0$ and $\omega_l=1$. Consequently, the hypothesis $(C_1)$, $(C_2)$ and $(C_3)$ imply that Proposition 2 can be applied repeatedly at least $l$ times. It follows that we have $a_m=\lambda_mT+\mu_m$, for $1\leq m < f(l+1)$. In order to have all the partial quotients of degree 1, up to infinity, we shall prove that  $\delta_m\omega_m\neq 0$ for all $m\geq l+1=f(1)$. For $m\geq f(1)$, we can write $m=f(n)+i$ where $n\geq 1$ and $0\leq i\leq r-1$. Therefore we want to prove that $\omega_{f(n)+i}\delta_{f(n)+i}\neq 0$ for $n\geq 1$ and  for $0\leq i\leq r-1$. To prove this, we shall show by induction that, for $n\geq 1$ and  for $0\leq i\leq r-1$, the following equalities hold :
 $$\delta_{f(n)}=a^{1-r}(\epsilon_{1,n}^{-1}\delta_n)^r \quad (D_{n,0})\quad \text{ and} \quad \delta_{f(n)+i}=(a/2)\lambda_{f(n)+i}^r\quad \text{ for} \quad 1\leq i\leq r-1, \eqno{(D_{n,i})}$$ 
together with
$$\omega_{f(n)+k}=\omega_n^r\quad (O_{n,k})\quad \text{ and }\quad \omega_{f(n)+i}=1\quad \text{for} \quad 0\leq i\leq r-1\quad \text{ and }\quad i\neq k.\eqno{(O_{n,i})}$$
  Note that if $\delta_j\omega_j\neq 0$ for $j<m=f(n)+i$, and $(D_{n,i})$ and  $(O_{n,i})$ hold, then we have $\delta_{m}\omega_{m}\neq 0$. 
\newline The proof of the equalities $(D_{n,i})$ and  $(O_{n,i})$ will follow by induction from $(DL_n)$, $(E_n)$, $(N_n)$, $(L_{n,i})$,$(M_{n,i})$ and $(X_{n,i})$.  
\newline First, we prove that $(D_{1,0})$ and  $(O_{1,0})$ hold. Using $(C_2)$, $(C_3)$, $(DL_{1})$ and $(L_{1,0})$, we have $\omega_l=1$, $\delta_l=-a\epsilon_1^r\epsilon_2^{-r}$, $\delta_1=a\lambda_1^r+\epsilon_2$ and $\lambda_{l+1}=\epsilon_1^{-1}\lambda_1^r$. Hence, from $(DL_{l+1})$, we get
$$ \delta_{f(1)}=\delta_{l+1}=a\lambda_{l+1}^r-a^{2-r}(\omega_l\delta_l)^{-1}=a\lambda_1^{r^2}\epsilon_1^{-r}+a^{1-r}\epsilon_2^r\epsilon_1^{-r}=a^{1-r}(\delta_1\epsilon_1^{-1})^r.$$
Besides, since $\omega_l=1$, we have $\pi_l=0$. Consequently, by $(N_{l+1})$, we get $\nu_{l+1}=0$. By $(M_{1,0})$, we have $\mu_{l+1}=0$. It follows that $\pi_{l+1}=0$ and $\omega_{f(1)}=\omega_{l+1}=1$.
\newline Let $n\geq 1$ and $0\leq i\leq r-1$. First, we shall prove that for $0\leq i\leq r-2$, $(D_{n,i})$ and  $(O_{n,i})$ imply $(D_{n,i+1})$ and $(O_{n,i+1})$. Secondly, we shall  prove that $(D_{n,r-1})$ and  $(O_{n,r-1})$ imply $(D_{n+1,0})$ and $(O_{n+1,0})$. The proof is divided into five cases, each comprising two parts.
\newline $\bullet$ {\it{Case 1: $i=0$. }} By $(D_{n,0})$, we have $\delta_{f(n)}=a^{-r+1}(\epsilon_{1,n}^{-1}\delta_n)^r$. Furthermore, by $(L_{n,1})$, we have $\lambda_{f(n)+1}=2(\delta_{n}\epsilon_{1,n}^{-1})^{-1}$.  Therefore,  with $\omega_{f(n)}=1$, from  $(DL_{f(n)+1})$, we get :
\begin{eqnarray*}
\delta_{f(n)+1}&=& a\lambda^r_{f(n)+1}-a^{2-r}(\omega_{f(n)}\delta_{f(n)})^{-1}\\
&=& a\lambda^r_{f(n)+1}-a^{2-r}a^{r-1}(\delta_n\epsilon_{1,n}^{-1})^{-r}\\
&=& a\lambda^r_{f(n)+1}-(a/2)\lambda^r_{f(n)+1}=(a/2)\lambda^r_{f(n)+1}.\\ 
\end{eqnarray*}
Hence $(D_{n,1})$ holds.
\newline By $(O_{n,0})$, we have $\omega_{f(n)}=1$ and, consequently, $\pi_{f(n)}=0$. By $(N_{f(n)+1})$, we obtain  $\nu_{f(n)+1}=0$. If $k>1$, $(M_{n,1})$ implies $\mu_{f(n)+1}=0$. Therefore, $\nu_{f(n)+1}=\mu_{f(n)+1}^r$ and $\omega_{f(n)+1}=1$. Thus $(O_{n,1})$ holds. If $k=1$ (i.e. $r=3$), the second part of Case 3 below must be applied.
\newline $\bullet$ {\it{Case 2: $i\neq 0,k-1,k,r-1$.}} By $(D_{n,i})$, we have $\delta_{f(n)+i}=(a/2)\lambda^r_{f(n)+i}$ and, by $(O_{n,i})$, $\omega_{f(n)+i}=1$. Furthermore, by $(X_{n,i})$, we have $\lambda_{f(n)+i}=4a^{-1}\lambda_{f(n)+i+1}^{-1}$.  Therefore, from  $(DL_{f(n)+i+1})$, we get :
\begin{eqnarray*}
\delta_{f(n)+i+1}&=& a\lambda^r_{f(n)+i+1}-a^{2-r}(\omega_{f(n)+i}\delta_{f(n)+i})^{-1}\\
&=& a\lambda^r_{f(n)+i+1}-a^{2-r}(2a^{-1}\lambda^{-r}_{f(n)+i})\\
&=& a\lambda^r_{f(n)+i+1}-2a^{1-r}(a^r/4)\lambda_{f(n)+i+1}^{r}=(a/2)\lambda^r_{f(n)+i+1}.\\
\end{eqnarray*} 
Hence $(D_{n,i+1})$ holds.
\newline By $(O_{n,i})$, we have $\omega_{f(n)+i}=1$ and consequently $\pi_{f(n)+i}=0$. By $(N_{f(n)+i+1})$, we obtain  $\nu_{f(n)+i+1}=0$. Since $i+1\neq k,k+1$, $(M_{n,i+1})$ implies $\mu_{f(n)+i+1}=0$. Therefore $\nu_{f(n)+i+1}=\mu_{f(n)+i+1}^r$ and $\omega_{f(n)+i+1}=1$. Thus $(O_{n,i+1})$ hold.
\newline $\bullet$ {\it{Case 3: $i=k-1$.}} If $k>1$, by the same arguments as in the first part of the previous case, since $i\neq k$, we see that $(D_{n,i})$ and  $(O_{n,i})$ imply $(D_{n,i+1})$. Hence $(D_{n,k})$ holds. If $k=1$ (i.e. $r=3$), the first part of Case 1 must be applied.
\newline By $(D_{n,k})$, we have $\delta_{f(n)+k}=(a/2)\lambda^r_{f(n)+k}$. By $(O_{n,k-1})$, we have $\omega_{f(n)+k-1}=1$ and $\pi_{f(n)+k-1}=0$. Hence, from $(N_{f(n)+k})$, we obtain $\nu_{f(n)+k}=0$. Therefore, using $(M_{n,k})$, we get :
\begin{eqnarray*}
\omega_{f(n)+k}&=& 1+a^{2-r}(\nu_{f(n)+k}-\mu^{r}_{f(n)+k})^2\delta_{f(n)+k}^{-2}\\
&=& 1+a^{2-r}4a^{-2}(\mu_{f(n)+k}\lambda^{-1}_{f(n)+k})^{2r}\\
&=& 1+4a^{-r}((-a)^{1-k}\pi_n/2)^{2r}\\
&=& 1+a^{2r-r^2}\pi_n^{2r}=(1+a^{2-r}\pi_n^2)^r=\omega_n^r.\\
\end{eqnarray*} 
Hence $(O_{n,k})$ holds.
\newline $\bullet$ {\it{Case 4: $i=k$.}} By $(D_{n,k})$, we have $\delta_{f(n)+k}=(a/2)\lambda^r_{f(n)+k}$ and, by $(O_{n,k})$, $\omega_{f(n)+k}=\omega_n^r$. Furthermore, by $(X_{n,k})$, we have $\omega_n\lambda_{f(n)+k}=4a^{-1}\lambda_{f(n)+k+1}^{-1}$.  Therefore, from  $(DL_{f(n)+k+1})$, we get :
\begin{eqnarray*}
\delta_{f(n)+k+1}&=& a\lambda^r_{f(n)+k+1}-a^{2-r}(\omega_{f(n)+k}\delta_{f(n)+k})^{-1}\\
&=& a\lambda^r_{f(n)+k+1}-a^{2-r}(2a^{-1}\omega_n^{-r}\lambda^{-r}_{f(n)+k})\\
&=& a\lambda^r_{f(n)+k+1}-2a^{1-r}(a^r/4)\lambda_{f(n)+k+1}^{r}=(a/2)\lambda^r_{f(n)+k+1}.\\
\end{eqnarray*} 
Hence $(D_{n,k+1})$ holds. 
\newline By $(O_{n,k-1})$, we have $\pi_{f(n)+k-1}=0$ and $\nu_{f(n)+k}=0$. By $(X_{n,k})$, we have $\omega_n\lambda_{f(n)+k}=4a^{-1}\lambda_{f(n)+k+1}^{-1}$.  By $(D_{n,k})$, we have $\delta_{f(n)+k}=(a/2)\lambda^r_{f(n)+k}$. Using  $(M_{n,k+1})$, from $(N_{f(n)+k+1})$, we get :
\begin{eqnarray*}
\nu_{f(n)+k+1}&=& a^{2-r}(\nu_{f(n)+k}-\mu^r_{f(n)+k})\omega_{f(n)+k}^{-1}\delta_{f(n)+k}^{-2}\\
&=& -a^{2-r}\mu^r_{f(n)+k}\omega_n^{-r}((a/2)\lambda^r_{f(n)+k})^{-2}\\
&=& -4a^{-r}\mu^r_{f(n)+k}(\omega_n\lambda_{f(n)+k})^{-r}\lambda_{f(n)+k}^{-r}\\
&=& -4a^{-r}\mu^r_{f(n)+k}(4a^{-1}\lambda_{f(n)+k+1}^{-1})^{-r}\lambda_{f(n)+k}^{-r}\\
&=& -(\mu_{f(n)+k}\lambda_{f(n)+k+1}\lambda_{f(n)+k}^{-1})^{r}=\mu^r_{f(n)+k+1}.\\
\end{eqnarray*}
Consequently, $\pi_{f(n)+k+1}=0$ and $\omega_{f(n)+k+1}=1$. Hence $(O_{n,k+1})$ holds.
\newline $\bullet$ {\it{Case 5: $i=r-1$.}}  Recall that $f(n+1)=f(n)+r$. By $(D_{n,r-1})$ and $(O_{n,r-1})$, we have $\delta_{f(n+1)-1}=\delta_{f(n)+r-1}=(a/2)\lambda_{f(n)+r-1}^r$ and $\omega_{f(n+1)-1}=\omega_{f(n)+r-1}=1$. By $(L_{n,r-1})$, we have $\lambda_{f(n)+r-1}=2a^{-1}\delta_{n}\epsilon_{1,n}^{-1}\omega_{n}^{(-1)^k}$ and also, by $(L_{n+1,0})$, $\lambda_{f(n+1)}=\lambda^{r}_{n+1}\epsilon_{1,n+1}^{-1}$. By $(E_{n+1})$, we also have $\epsilon_{1,n}^{-1}\omega_{n}^{(-1)^k}=a^{r-1}\epsilon_{1,n+1}\omega_{n}$. Moreover, by $(DL_{n+1})$, we have $(\omega_{n}\delta_{n})^{-1}=a^{r-2}(a\lambda_{n+1}^r-\delta_{n+1})$. Therefore, from  $(DL_{f(n+1)})$, we get :
\begin{eqnarray*}
\delta_{f(n+1)}&=& a\lambda^r_{f(n+1)}-a^{2-r}(\omega_{f(n+1)-1}\delta_{f(n+1)-1})^{-1}\\
&=& a\lambda^r_{f(n+1)}-a^{2-r}(2a^{-1}\lambda_{f(n)+r-1}^{-r})\\
&=& a\lambda^r_{f(n+1)}-2a^{1-r}(2a^{-1}\delta_{n}\epsilon_{1,n}^{-1}\omega_{n}^{(-1)^k})^{-r}\\
&=& a\lambda^r_{f(n+1)}-a(a^{r-1}\epsilon_{1,n+1}\omega_{n}\delta_{n})^{-r}\\
&=& a\lambda^{r^2}_{n+1}\epsilon_{1,n+1}^{-r}-a^{-r^2+r+1}\epsilon_{1,n+1}^{-r}(a^{r-2}(a\lambda_{n+1}^r-\delta_{n+1}))^{r}\\
&=& \epsilon_{1,n+1}^{-r}(a\lambda^{r^2}_{n+1}-a^{1-r}(a^r\lambda_{n+1}^{r^2}-\delta_{n+1}^r))= a^{1-r}(\epsilon_{1,n+1}^{-1}\delta_{n+1})^r.\\
\end{eqnarray*}
Consequently $(D_{n+1,0})$ holds.
\newline By $(O_{n,r-1})$, we have $\pi_{f(n)+r-1}=0$. Thus, from $(N_{f(n+1)})$, we get $\nu_{f(n+1)}=0$. By $(M_{n+1,0})$, we have $\mu_{f(n+1)}=0$. Therefore, $\pi_{f(n+1)}=0$ and $\omega_{f(n+1)}=1$. Consequently, $(O_{n+1,0})$ holds. 
\par Thus we have proved that $\delta_m\omega_m\neq 0$, for $m\geq 1$ and this implies that $\alpha \in \mathcal{E}^*(r,l,a,q)$. 
\newline Now we turn to the description of the sequence of partial quotients $a_n=\lambda_nT+\mu_n$. The first $l$ values of both sequences $(\lambda_n)_{n\geq 1}$ and $(\mu_n)_{n\geq 1}$ are given, as well as $\delta_i$, $\nu_i$ and $\omega_i$ for $1\leq i\leq l$. We recall that, for $n\geq 1$, we have $g(n)=nr+l-k=f(n)+k$. According to $(O_{n,k})$, we have $\omega_{g(n)}=\omega_n^r$. According to $(O_{n,i})$ for $i\neq k$, if $m>l$ and $m\neq g(n)$, we have $\omega_m=1$. For $i\in I$, we have $\omega_i=1$ if and only if $i\notin I^*$. Consequently, by iteration, we obtain 
$$\omega_n=1 \quad \text{if}\quad n\notin G\cup I^*\quad \text{and}\quad \omega_n=\omega_i^{r^m} \quad \text{if}\quad n=g^m(i)\quad \text{ for}\quad m\geq 1\quad \text{and}\quad i\in I^*.\eqno{(35)}$$
We start by the description of the sequence $(\mu_n)_{n\geq l+1}$. For $n\geq l+1$, if $n\notin G$, we have $\omega_n=1$ and therefore $\pi_n=0$. This implies, according to $(M_{n,i})$ for $0\leq i\leq r-1$, that $\mu_n=0$ if $n\notin  G\cup (G+1)$. Let $i\in I^*$ and $m\geq 1$, since $\nu_{n+1}=a^{2-r}\pi_n(\omega_n\delta_n)^{-1}$ and $\pi_{g^m(i)-1}=0$, we get $\nu_{g^m(i)}=0$, and we have
$$\pi_{g^m(i)}=-\mu^r_{g^m(i)}\delta^{-1}_{g^m(i)} \quad \text{ for}\quad m\geq 1\quad \text{and}\quad i\in I^*.\eqno{(36)}$$ 
We set $n'=g^{m-1}(i)$, then $g^m(i)=f(n')+k$. Applying $(M_{n',k})$ and $(D_{n',k})$, we have
$$ (\mu_{g^m(i)}\lambda^{-1}_{g^m(i)})^r=(-a)^{(1-k)r}\pi_{n'}^r/2 \quad \text{and} \quad \delta^{-1}_{g^m(i)} =2a^{-1}\lambda^{-r}_{g^m(i)}.\eqno{(37)}$$
From $(36)$ and $(37)$, we get
$$\pi_{g^m(i)}=-2a^{-1}(\mu_{g^m(i)}\lambda^{-1}_{g^m(i)})^r=(-a)^{(1-k)r-1}\pi_{n'}^r.\eqno{(38)}$$
We set $A=(-a)^{(1-k)r-1}$. Then $(38)$ becomes $\pi_{g^m(i)}=A\pi_{g^{m-1}(i)}^r$. By iteration, we get 
$$\pi_{g^m(i)}=A^{u_m}\pi_i^{r^m}\quad \text{ for}\quad m\geq 1\quad \text{and}\quad i\in I^*,\quad \text{ where}\quad u_m=(r^m-1)/(r-1).\eqno{(39)}$$  
Hence, if $n=g^{m}(i)$ for $i\in I^*$ and $m\geq 1$, by $(M_{n',k})$ and $(39)$, we have
$$\mu_n\lambda_n^{-1}=(-a)^{1-k}\pi_{g^{m-1}(i)}/2=(-a)^{1-k}A^{u_{m-1}}\pi_i^{r^{m-1}}/2=(-a)^{v_m} ((\nu_i-\mu_i^r)\delta_i^{-1})^{r^{m-1}}/2,\eqno{(40)}$$
where $v_m=(1-k)+((1-k)r-1)u_{m-1}$. It is elementary to check that $v_m=(r^{m-1}(2-r)+1)/2$. Recalling that, by $(M_{n',k+1})$, we also have $\mu_n\lambda_n^{-1}=-\mu_{n+1}\lambda_{n+1}^{-1}$, with $(40)$ we have completed the description of the sequence $(\mu_n)_{n\geq 1}$. 
\newline Finally we turn to the description of the sequence $(\lambda_n)_{n\geq l+1}$. We recall that, according to $(L_{1,0})$, we have $\lambda_{l+1}=\epsilon_1^{-1}\lambda_1^r$. Hence, this description will follow from computing $C(n)=\lambda_{n-1}\lambda_n$ for $n>l+1$. If $n=g^{m}(i)+1$ for $i\in I^*$ and $m\geq 1$, then $n-1=f(n')+k$, where $n'=g^{m-1}(i)$. Consequently, applying $(X_{n',k})$ and $(35)$, we have
$$\lambda_n\lambda_{n-1}=\lambda_{f(n')+k+1}\lambda_{f(n')+k}=4a^{-1}\omega_{n'}^{-1}=4a^{-1}\omega_i^{-r^{m-1}}=4a^{-1}(1+a^{2-r}(\nu_i-\mu_i^r)^2\delta_i^{-2})^{-r^{m-1}} .\eqno{(41)}$$
According to $(L_{n,0})$ and $(D_{n,0})$, we have $\lambda_{f(n)}=\epsilon_{1,n}^{-1}\lambda_n^r$ and $\delta_{f(n)}=a^{1-r}(\epsilon_{1,n}^{-1}\delta_n)^r$. Hence, we obtain directly
$$\delta_{f(n)}^{-1}\lambda_{f(n)}^{r}=a^{r-1}(\delta_n^{-1}\lambda_n^{r})^r .\eqno{(42)}$$
For $i\in I$ and $m\geq 1$, from $(42)$, we get $\delta_{f^m(i)}^{-1}\lambda_{f^m(i)}^{r}=a^{r-1}(\delta_{f^{m-1}(i)}^{-1}\lambda_{f^{m-1}(i)}^{r})^r$. Consequently, by iteration, we get 
$$\delta_n^{-1}\lambda_n^{r}=a^{r^m-1}(\delta_i^{-1}\lambda_i^{r})^{r^m} \quad \text{if}\quad n=f^m(i)\quad \text{ for}\quad m\geq 1\quad \text{and}\quad i\in I.\eqno{(43)}$$
We set $n'=f^{m-1}(i)$. By $(L_{n',0})$ and $(L_{n',1})$, since $\lambda_{f(n')+1}=2(\epsilon_{1,n'}^{-1}\delta_{n'})^{-1}$, from $(43)$, we obtain
$$\lambda_{f(n')+1}\lambda_{f(n')}=2\delta_{n'}^{-1}\lambda_{n'}^{r} =2a^{r^{m-1}-1}(\delta_i^{-1}\lambda_i^{r})^{r^{m-1}}.\eqno{(44)}$$
 Hence, by $(44)$, we have
$$\lambda_n\lambda_{n-1}=2a^{-1}(a\lambda_i^{r}\delta_i^{-1})^{r^{m-1}} \quad \text{if}\quad n=f^m(i)+1\quad \text{ for}\quad m\geq 1\quad \text{and}\quad i\in I.\eqno{(45)}$$
We set $n=f^{m}(i)$ for $m\geq 1$ and $i\in I$, with $n\neq f(1)$.  We have $n=f(n')$ with $n'>1$. Consequently, by $(D_{n'-1,r-1})$, we get
$$\delta_{n-1}=\delta_{f(n')-1}=\delta_{f(n'-1)+r-1}=(a/2)\lambda_{f(n'-1)+r-1}=(a/2)\lambda_{n-1}^{r}.\eqno{(46)}$$
As $n>l+1$, then $n-1\geq l+1$ and $n-1\notin G$. Therefore, by $(35)$, we have $\omega_{n-1}=1$. Consequently, by $(L_{n,0})$ and by $(E_n)$, since  $\omega_{n-1}=1$, we have 
$$ \lambda_n^r=\epsilon_{1,n} \lambda_{f(n)} =a^{1-r}\epsilon_{1,n-1}^{-1}\lambda_{f(n)} . \eqno{(47)}$$
Combining $(46)$ and $(47)$, we obtain
$$(a\lambda_{n-1}\lambda_{n})^r=a^r(2a^{-1}\delta_{n-1})a^{1-r}\epsilon_{1,n-1}^{-1}\lambda_{f(n)} =2\delta_{n-1}\epsilon_{1,n-1}^{-1}\lambda_{f(n)}.\eqno{(48)}$$
Also, by $(L_{n-1,r-1})$, since  $\omega_{n-1}=1$, we have
$$\lambda_{f(n)-1}=\lambda_{f(n-1)+r-1}=2a^{-1}\delta_{n-1}\epsilon_{1,n-1}^{-1}.\eqno{(49)}$$
Combining $(48)$ and $(49)$, we obtain
$$(a\lambda_{n-1}\lambda_{n})^r=a\lambda_{f(n)-1}\lambda_{f(n)}.\eqno{(50)}$$
Hence, by iteration from $(50)$, we get
$$a\lambda_{f^m(i)-1}\lambda_{f^m(i)}=(a\lambda_{f(i)-1}\lambda_{f(i)})^{r^{m-1}}\quad \text{ for}\quad m\geq 1\quad \text{and}\quad i\neq 1\in I\eqno{(51)}$$
and also
$$a\lambda_{f^m(1)-1}\lambda_{f^m(1)}=(a\lambda_{f^2(1)-1}\lambda_{f^2(1)})^{r^{m-2}}\quad \text{ for}\quad m\geq 2.\eqno{(52)}$$
Now, let $n\geq 2$, as above we have   
$$\lambda_{f(n)}=\epsilon_{1,n}^{-1}\lambda_n^r\quad \text{and}\quad\lambda_{f(n)-1}=2a^{-1}\delta_{n-1}\epsilon_{1,n-1}^{-1}\omega_{n-1}^{(-1)^k}.\eqno{(53)}$$
By $(E_{n})$, we have 
$$\omega_{n-1}^{(-1)^k}=a^{r-1}\omega_{n-1}\epsilon_{1,n-1}\epsilon_{1,n}.\eqno{(54)}$$
From $(53)$ and $(54)$, we obtain
 $$\lambda_{f(n)-1}\lambda_{f(n)}=2a^{r-2}\delta_{n-1}\omega_{n-1}\lambda_n^r.\eqno{(55)}$$
Recalling that, by $(DL_n)$, we have $a^{2-r}(\delta_{n-1}\omega_{n-1})^{-1}=a\lambda_n^r-\delta_n$, we obtain
$$\lambda_{f(n)-1}\lambda_{f(n)}=2\lambda_n^r(a\lambda_n^r-\delta_n)^{-1}=2a^{-1}(1-a^{-1}\lambda_n^{-r}\delta_n)^{-1}.\eqno{(56)}$$
Combining $(51)$ and $(56)$, we get
$$\lambda_{f^m(i)-1}\lambda_{f^m(i)}=2a^{-1}(1-a^{-1}\lambda_i^{-r}\delta_i)^{-r^{m-1}}\quad \text{ for}\quad m\geq 1\quad \text{and}\quad i\neq 1\in I.\eqno{(57)}$$
Besides, according to $(56)$, we have
$$a\lambda_{f^2(1)-1}\lambda_{f^2(1)}=2(1-a^{-1}\lambda_{f(1)}^{-r}\delta_{f(1)})^{-1}.\eqno{(58)}$$
From $(42)$, we get
$$a^{-1}\lambda_{f(1)}^{-r}\delta_{f(1)}=(a^{-1}\lambda_1^{-r}\delta_1)^r.\eqno{(59)}$$
Hence, by $(58)$ and $(59)$, we have
$$a\lambda_{f^2(1)-1}\lambda_{f^2(1)}=2(1-a^{-1}\lambda_1^{-r}\delta_1)^{-r}.\eqno{(60)}$$
Combining $(52)$ and $(60)$, we get
$$\lambda_{f^m(1)-1}\lambda_{f^m(1)}=2a^{-1}(1-a^{-1}\lambda_1^{-r}\delta_1)^{-r^{m-1}}\quad \text{ for}\quad m\geq 2.\eqno{(61)}$$
In summary, according to $(41)$, $(45)$, $(57)$ and $(61)$, we have obtained the values stated in the theorem for $C(n)$ if $n\in F\cup (F+1)\cup (G+1)$. Let us turn to the last case $n\notin F\cup (F+1)\cup (G+1)$. We have $n>l+1$, consequently there exist $n_1\geq 1$ and $i=0,\dots,r-1$ such that $n-1=f(n_1)+i$ and $n-1\notin (F-1)\cup F\cup G$. There will be four cases according to the value of $i$.
\newline $\bullet$ {\it{Case 1: $i\neq 0,k,r-1$.}} According to $(X_{n_1,i})$, we have
$$\lambda_{n-1}\lambda_{n}=\lambda_{f(n_1)+i}\lambda_{f(n_1)+i+1}=4a^{-1}.\eqno{(62)}$$
 $\bullet$ {\it{Case 2: $i=k$.}} Here we have $n-1=g(n_1)$ and , since $n-1\notin G$, we have $n_1\notin G\cup I^*$. Therefore, by $(35)$, we have $\omega_{n_1}=1$. Consequently, according to $(X_{n_1,k})$, we have 
$$\lambda_{n-1}\lambda_{n}=\lambda_{f(n_1)+k}\lambda_{f(n_1)+k+1}=4a^{-1}\omega_{n_1}^{-1}=4a^{-1}.\eqno{(63)}$$
For the last two cases, $i=0$ or $i=r-1$, we need the following. If $n\notin F\cup I$ then there are three integers $m\geq 0$, $n_2\geq 1$ and $j=1,\dots,r-1$ such that $n=f^m(f(n_2)+j)$. Set $n_3=f(n_2)+j$. By $(42)$ and by iteration, since $n=f^m(n_3)$, we can write
$$\delta_{n}^{-1}\lambda_{n}^r=a^{r^m-1}(\delta_{n_3}^{-1}\lambda_{n_3}^r)^{r^m}.\eqno{(64)}$$
By $(D_{n_2,j})$, with $j\neq 0$, we get $\delta_{n_3}^{-1}\lambda_{n_3}^r=2a^{-1}$. Consequently, by $(64)$, the previous argument implies : 
$$n \notin F\cup I\quad \Rightarrow \quad \delta_{n}^{-1}\lambda_{n}^r=2a^{-1}.\eqno{(65)}$$
$\bullet$ {\it{Case 3: $i=0$.}} Here we have $n-1=f(n_1)$. Hence, according to $(44)$, we have 
$$\lambda_{n-1}\lambda_{n}=\lambda_{f(n_1)}\lambda_{f(n_1)+1}=2\delta_{n_1}^{-1}\lambda_{n_1}^r.\eqno{(66)}$$
Since $n-1\notin F$ and $n-1=f(n_1)$, we have $n_1\notin F\cup I$. Therefore, by $(65)$, we have  $\delta_{n_1}^{-1}\lambda_{n_1}^r=2a^{-1}$. Consequently, $(66)$ becomes $\lambda_{n-1}\lambda_{n}=4a^{-1}$. 
\newline $\bullet$ {\it{Case 4: $i=r-1$.}}  Here we have $n=f(n_1)+r=f(n_1+1)$. According to $(56)$, we have
$$\lambda_{n-1}\lambda_{n}=\lambda_{f(n_1+1)-1}\lambda_{f(n_1+1)}=2a^{-1}(1-a^{-1}\lambda_{n_1+1}^{-r}\delta_{n_1+1})^{-1}.\eqno{(67)}$$
Since $n\notin F$ and $n=f(n_1+1)$, we have $n_1+1\notin F\cup I$. Therefore, by $(65)$ we have  $\lambda_{n_1+1}^{-r}\delta_{n_1+1}=a/2$. Hence, from $(67)$, we get $$\lambda_{n-1}\lambda_{n}=2a^{-1}(1-a^{-1}\lambda_{n_1+1}^{-r}\delta_{n_1+1})^{-1}=2a^{-1}(1-a^{-1}a/2)^{-1}=4a^{-1}.\eqno{(68)}$$ 
In summary, according to $(62)$, $(63)$, $(66)$ and $(68)$, if $n\notin F\cup (F+1)\cup (G+1)$, we have obtained $C(n)=4a^{-1}$. Hence the description of the sequence $(\lambda_n)_{n\geq 1}$ is over.
\newline So the proof of the theorem is complete. 

\vskip 0.5 cm
\noindent{\bf{5. Last comments}}
\par We want to come back to the statement of the theorem presented in this note. Starting from $\alpha \in \mathcal{E}(r,l,a,q)$, we have proved that the three conditions $(C_1)$, $(C_2)$ and $(C_3)$ are sufficient to have $\alpha \in \mathcal{E}^*(r,l,a,q)$. However, condition $(C_1)$ is particular and clearly necessary. Indeed, by repeated application of Proposition 2, it allows to have the first (up to the rank $l+r(l-1)$) partial quotients of degree 1. While conditions $(C_2)$ and $(C_3)$ are useful to keep this repetition of Proposition 2 up to infinity and to obtain a sequence of partial quotients having a relatively simple pattern. Then it is natural to ask whether conditions $(C_2)$ and $(C_3)$ are also necessary to have $\alpha \in \mathcal{E}^*(r,l,a,q)$. We know that the subset $\mathcal{E}(r,l,a,q)$ is finite, consequently all the elements can be tested by computer. We have done so for $r=q=p$ and for small values of $p$ and $l$. We have observed that, for $\alpha \in \mathcal{E}(p,l,a,p)$, if $(C_1)$, $(C_2)$ or $(C_3)$ is not satisfied then $\alpha \notin \mathcal{E^*}(p,l,a,p)$. Consequently we conjecture that this set of conditions is not only sufficient but also necessary when the finite base field is prime. However, this is not generally so if the base field is not prime. Indeed, it was a surprise to discover in $\mathcal{E}_0(r,l,a,q)$, with $q>p$, certain continued fractions belonging to  $\mathcal{E}^*_0(r,l,a,q)$ but for which the condition $(C_2)$ of the theorem stated here is not satisfied. This phenomenon has been explained in [7] and an example in $\mathcal{E}^*_0(3,1,2,27)$ has been given there [7, p. 258]. Note that for this type of example the sequence of partial quotients has a much more complex pattern than the ones presented in this note.
\par With the conditions of our theorem, let us give a short and global description of the structure of this sequence of partial quotients.  We can write this infinite continued fraction expansion $\alpha=[a_1,a_2,\dots,a_n,\dots]$ in the following way :
$$\alpha=[W_{\infty}]=[W_0,W_1,\dots,W_m,\dots],$$
where, for $m\geq 0$, we define the finite word  $$W_m=a_{f^m(1)},a_{f^m(1)+1},\dots,a_{f^{m+1}(1)-1}.$$
Note that we have $W_0=a_1,\dots,a_l$, which plays the role of the "basis" of the expansion. We also have $W_1=a_{l+1},\dots,a_{lr+l}$ and so on. We can check that $\vert W_m \vert =lr^m$ for $m\geq 0$. Define, for $n\geq 1$, the following word of lenght $r$ :
 $$W^f(a_n)=a_{f(n)},a_{f(n)+1},\dots, a_{f(n)+r-1}.$$
 Then with these notations, we see that $W_{m+1}$ is built upon $W_m$, by concatenation, in the following way
$$W_{m+1}=W^f(a_{f^m(1)}),W^f(a_{f^m(1)+1}),\dots,W^f(a_{f^{m+1}(1)-1}).$$ 
We see that the pattern of these sequences is somehow very regular. Let us come back to the original examples, due to Mills and Robbins [12, p. 400]. They belong to $\mathcal{E}^*_0(p,2,4,p)$ for all $p>3$ and they are defined by $(\lambda_1,\lambda_2,\epsilon_1,\epsilon_2)=(\lambda_1,-\lambda_1(1+2\lambda_1)^{-1},1,2)$, with $\lambda_1\in \F_p^*$ and $\lambda_1\neq -1/2$. Here, we have $G=\emptyset$, $\omega_n=1$ and $\mu_n=0$ for $n\geq 1$. The conditions of the theorem are satisfied and the description of the sequence $(\lambda_n)_{n\geq 1}$ follows. In 1988, shortly after the publication of these examples, J-P. Allouche [1] could show the regularity of this sequence $(\lambda_n)_{n\geq 1}$, by proving that it is $p$-automatic. It should probably be possible to extend this type of result to all the sequences presented in this note.

\vskip 0.5 cm
\noindent{\bf{Acknowledgment}}
\par During the winter 2011-2012, the first author kept up a correspondence with
Alina Firicel who had a postdoc position at Grenoble university. The
matter presented in this note was discussed and a first version of
Proposition 1 was prepared. We wish to thank her for this
contribution and, more generally, for her interest in this work. Part of the work was done while
the second author visited the Institut de Math\'ematiques de Bordeaux,
and he would like to thank his colleagues for their generous hospitality.
He would also like to thank the National Natural Science Foundation
of China (Grants no. 10990012 and 11371210) and the Morningside Center of
Mathematics (CAS) for partial financial support.

\begin{tabular}{ll}
 Alain LASJAUNIAS\\
Institut de Math\'ematiques de Bordeaux  \\
CNRS-UMR 5251\\
Talence 33405 \\
France \\
E-mail: Alain.Lasjaunias@math.u-bordeaux1.fr\\
\\
Jia-Yan YAO\\
Department of Mathematics\\
Tsinghua University\\
Beijing 100084\\
People's Republic of China\\
E-mail: jyyao@math.tsinghua.edu.cn
\end{tabular}%


\vskip 1 cm
\begin{thebibliography}{7}  

\addcontentsline{toc}{section}{References}
\bibitem[1]{} J.-P. Allouche, \emph{Sur le d\'eveloppement en fraction continue de certaines s\'eries
formelles}, C. R. Acad. Sci. Paris {\bf 307} (1988), 631-633.

\bibitem[2]{} L. Baum and M. Sweet,\emph{ Continued fractions of algebraic power series in characteristic 2}, Annals of Mathematics {\bf 103}
  (1976), 593-610.

\bibitem[3]{} L. Baum and M. Sweet,\emph{ Badly approximable power series in characteristic 2}, Annals of Mathematics {\bf 105}
  (1977), 573-580. 

\bibitem[4]{} A. Bluher and A. Lasjaunias,\emph{ Hyperquadratic power
    series of degree four}, Acta Arithmetica {\bf 124}
  (2006), 257-268.

\bibitem[5]{} D. Gomez and A. Lasjaunias,\emph{ Hyperquadratic power
    series in $\F_3((T^{-1}))$ with partial quotients of degree 1}, The Ramanujan Journal {\bf 33.2} (2014), 219-226.  

\bibitem[6]{} A. Lasjaunias,\emph{ Continued fractions for
    hyperquadratic power series over a finite field},  Finite Fields and their
  Applications {\bf 14} (2008), 329--350. 

 
\bibitem[7]{} A. Lasjaunias, \emph{ Algebraic continued fractions in
    $\F_q((T^{-1}))$ and recurrent sequences in $\F_q$}, Acta
  Arithmetica {\bf 133.3} (2008), 251-265.  

\bibitem[8]{} A. Lasjaunias,\emph{ Quartic power series in
    $\F_3((T^{-1}))$ with bounded
    partial quotients}, Acta Arithmetica {\bf 95}
  (2000), 49-59.

\bibitem[9]{} A. Lasjaunias and B. de Mathan,\emph{ Differential equations and diophantine approximation in positive characteristic}, Monatshefte fur Mathematik {\bf 128}
  (1999), 1-6.

\bibitem[10]{} A. Lasjaunias and J-J. Ruch,\emph{ Flat power series over
    a finite field}, 
Journal of Number Theory {\bf 95} (2002), 268--288.  

\bibitem[11]{} M. Mend\`es France, \emph{ Sur les fractions continues limit\'ees}, Acta
  Arithmetica {\bf 23} (1973), 207-215.  


\bibitem[12]{} W. Mills and D. Robbins,\emph{ Continued fractions for certain algebraic power series}, 
Journal of Number Theory {\bf 23} (1986), 388--404.

\bibitem[13]{} O. Perron, \emph{Die Lehre von den Kettenbr\"{u}chen}, 2nd
  ed., Chelsea Publishing Company, New York, 1950.

\bibitem[14]{} W. Schmidt, \emph{ On continued fractions and diophantine
    approximation in power series fields}, Acta
  Arithmetica {\bf 95} (2000), 139-166.  
\bibitem[15]{} D. Thakur, \emph{ Function Field Arithmetic}, World Scientific, 2004. 


\end{thebibliography}
\end{document}